\documentclass[10pt]{amsart}
\usepackage{latexsym, amsmath}

 \renewcommand{\epsilon}{\varepsilon}
 \newcommand{\newsection}[1]
 {\subsection{#1}\setcounter{theorem}{0} \setcounter{equation}{0}
 \par\noindent}

 \newtheorem{theorem}{Theorem}
 
 \newtheorem{lemma}[theorem]{Lemma}
 \newtheorem{corr}[theorem]{Corollary}
 
 \newtheorem{proposition}[theorem]{Proposition}
 \newtheorem{deff}[theorem]{Definition}
 
 \newcommand{\bth}{\begin{theorem}}
 \newcommand{\ble}{\begin{lemma}}
 \newcommand{\bcor}{\begin{corr}}
 \newcommand{\bdeff}{\begin{deff}}
 \newcommand{\bprop}{\begin{proposition}}
 \newcommand{\eth}{\end{theorem}}
 \newcommand{\ele}{\end{lemma}}
 \newcommand{\ecor}{\end{corr}}
 \newcommand{\edeff}{\end{deff}}
 
 \newcommand{\eprop}{\end{proposition}}
  
  \newcommand{\Rn}{{\Bbb{R}^n}}

 \newcommand{\supp}{\text{supp }}
 \renewcommand{\Pi}{\varPi}

 \renewcommand{\epsilon}{\varepsilon}

 \newcommand{\Rplus}{{\Bbb R}_+}
 \newcommand{\R}{{\Bbb R}}

\renewcommand{\square}{\Box}

\begin{document}

 \title[Global existence for semilinear wave equations]
 {Global existence for semilinear wave equations exterior to nontrapping obstacles}
\author{Jason L. Metcalfe}
\thanks{The author would like to thank C. Sogge and A. Stewart for numerous helpful discussions during this study.}
\address{Department of Mathematics, The Johns Hopkins University,
Baltimore, MD 21218}
\subjclass[2000]{Primary 35L05}
\date{October 17, 2002}
\begin{abstract}
In this paper, we prove the existence of global small amplitude
solutions to semilinear wave equations with quadratic
nonlinearities exterior to a nontrapping obstacle.  This
generalizes the work of Hayashi in a domain exterior to a ball and
of Shibata and Tsutsumi in spatial dimensions $n\ge 6$.
\end{abstract}

 \maketitle

\newsection{Introduction}

The goal of this paper is to prove global existence of solutions
to semilinear wave equations with quadratic nonlinearities
exterior to a nontrapping obstacle.  This extends results that
were previously known in Minkowski space (see, e.g., \cite{S}).
\par
More precisely, let $\mathcal{K}\subset \Rn$ be a smooth, compact,
nontrapping obstacle, and set $\Omega=\Rn\backslash \mathcal{K}$.
We shall consider solutions to
\begin{equation}\label{main}
\begin{cases}
\Box u=\partial_t^2u-\Delta_x u=Q(u'), \quad (t,x)\in \Rplus\times \Omega\\
u(t,x)=0\quad \text{for } x\in\partial\Omega\\
u(0,\cdot)=f,\quad \partial_tu(0,\cdot)=g.
\end{cases}
\end{equation}
Here $Q$ is a constant coefficient quadratic form in
$u'=(\partial_t u, \nabla_x u)$.
\par
In order to solve \eqref{main}, we need to assume that the data
satisfies certain compatibility conditions.  To describe these
briefly, let $J_ku=\{\partial^\alpha_x u : 0\le|\alpha|\le k\}$.
Then, for fixed $m$ and $u$ a formal $H^m$ solution of
\eqref{main}, we can write $\partial_t^k u(0,\cdot)=\psi_k(J_k f,
J_{k-1}g)$, $0\le k\le m$, for compatibility functions $\psi_k$
which depend on the nonlinear term $Q$, $J_k f$, and $J_{k-1}g$.
The compatibility condition for \eqref{main} with $(f,g)\in
H^m\times H^{m-1}$ requires that $\psi_k$ vanish on
$\partial\Omega$ when $0\le k\le m-1$.  Additionally, one says
that $(f,g)\in C^\infty$ satisfy the compatibility condition to
infinite order if this condition holds for all $m$.  For more
detail on compatibility conditions, see e.g., \cite{KSS1}.
\par
In proving the main theorem, we will need to use the invariance of
the wave equation under translations and spatial rotations.  Let
$Z=\{\partial_t, \partial_k, x_i\partial_j - x_j\partial_i\}$,
$1\le k\le n$, $1\le i<j\le n$. The vector fields of $Z$
essentially preserve the Dirichlet boundary conditions.  The
Lorentz boosts $\Omega_{0k}=t\partial_k+x_k\partial_t$ do not have
this property, and thus, do not seem appropriate for use in
obstacle problems.
\par
We are now ready to state the main result.

\begin{theorem}\label{GEexterior}  Suppose $n\ge 4$.
Let $\Omega$ be a domain in $\Rn$ exterior to a smooth compact
nontrapping obstacle and assume that $Q(u')$ is as above.  Assume
that $(f,g)\in C^\infty (\Omega)$ satisfies the compatibility
conditions to infinite order. If $\varepsilon>0$ is small and
\begin{equation}\label{GEexteriorAssumption}
\sum_{|\alpha|\le n+2}\|Z^\alpha
f\|_{L^2(\Omega)}+\sum_{|\alpha|\le n+1}\|Z^\alpha
g\|_{L^2(\Omega)}\le \varepsilon,
\end{equation}
then \eqref{main} has a unique global solution $u\in C^\infty
(\Rplus\times\Omega)$.
\end{theorem}

Theorem \ref{GEexterior} is an extension of previous results of
Shibata and Tsutsumi \cite{ST} and Hayashi \cite{Hayashi}. Shibata
and Tsutsumi were able to prove Theorem \ref{GEexterior} when
$n\ge 6$.  Their method, however, breaks down for the cases
$n=4,5$ and requires that the nonlinearity be cubic.  Hayashi was
able to extend the result to $n\ge 4$ when the obstacle
$\mathcal{K}$ is a ball.  Here, we extend the result to all
domains exterior to a nontrapping obstacle.  The techniques that
we use are similar to those used by Keel, Smith, and Sogge in
\cite{KSS2, KSS3} to show that well-known almost global existence
results in $n=3$ for solutions of the Minkowski wave equation
extend to exterior domains.
\par
We will begin by showing that one can prove global existence in
Minkowski space using these techniques.  Specifically, we will
consider
\begin{equation}\label{mainMinkowski}
\begin{cases}
\Box u = Q(u'),\quad (t,x)\in \Rplus\times\Rn\\
u(0,\cdot)=f,\quad \partial_t u(0,\cdot)=g.
\end{cases}
\end{equation}
In this case, we shall prove
\begin{theorem}\label{GEminkowski}
Let $n\ge 4$.  Assume that $Q(u')$ is as above and that $(f,g)\in
C^\infty(\Rn)$.  If $\varepsilon>0$ is small and
\begin{equation}\label{GEminkowskiAssumption}
\sum_{|\alpha|\le n+2}\|Z^\alpha f\|_{L^2(\Rn)}+\sum_{|\alpha|\le
n+1}\|Z^\alpha g\|_{L^2(\Rn)}\le \varepsilon,
\end{equation}
then \eqref{mainMinkowski} has a unique global solution $u\in
C^{\infty}(\Rplus\times\Rn)$.
\end{theorem}

This paper is organized as follows.  In the following section, we
will prove the main estimates in Minkowski space that we will
need.  In Section 3, we will prove Theorem \ref{GEminkowski}.  In
Section 4, we will extend the results of Section 2 and collect the
main estimates that will be needed to show global existence in the
exterior domain.  Finally, in Section 5, we will prove Theorem
\ref{GEexterior}.
\par
In a future paper, we hope to extend these results via similar
techniques to handle quasilinear equations.


\newsection{Main estimates in free space.}
We will begin by proving a few results related to the standard
energy inequality.  That is, if $v$ is a solution to the wave
equation $\square v=G$ in $\Rplus \times \Rn$, we have
\begin{equation}\label{stenergy}
\|v'(t,\cdot)\|_{L^2(\Rn)}\le
C\|v'(0,\cdot)\|_{L^2(\Rn)}+C\int_0^t
\|G(s,\cdot)\|_{L^2(\Rn)}\:ds.
\end{equation}

The first of these is a weighted energy estimate.  This result
follows from a simple modification of the arguments in Hormander
\cite{H} (Lemma 6.3.5).  A detailed proof can be found in
\cite{Metcalfe}.

\begin{lemma}\label{weightedEnergy}
Suppose that $n\geq 4$. Let $v(t,x)$ be a solution to the
homogeneous Minkowski wave equation $\Box v=0$ with initial data
$f,g\in C_C^{\infty}(\Rn)$ supported in $\{|x|\leq R\}$.  Then,
the following estimate holds
$$\int (t-|x|+2)^2\left(\bigl|\nabla_x v(t,x)\bigr|^2 + (\partial_t v(t,x))^2\right)\:dx
\leq C_R \left(\int \bigl|\nabla f\bigr|^2 + |g|^2\:dx\right).$$
\end{lemma}

We will also need $L^2_tL^2_x$ estimates.  The first of these is
an $L^2_tL^2_x\to L^2_tL^2_x$ estimate when the forcing term is
supported in a ball of fixed radius for any time $t$.

\begin{proposition}\label{prop2}
Suppose $n\ge 4$.  Let $v$ be a solution of the wave equation
$\square v=G$ in $\Rplus\times\Rn$ with vanishing Cauchy data.
Suppose, also, that $G(s,x)=0$ when $|x|>2$. Then,
$$\sum_{|\alpha|\le N}\|(1+r)^{-(n-1)/4}Z^{\alpha}v'(s,x)\|_{L^2_{s,x}([0,t]\times\Rn)}
\le C\sum_{|\alpha|\le N}\|Z^\alpha
G\|_{L^2_{s,x}([0,t]\times\Rn)}
$$
for a uniform constant $C$.
\end{proposition}

The proof of this proposition will be based on the following
lemma.
\begin{lemma}\label{prop2lemma}
Let $v$ be a solution of the wave equation $\square v=0$ in
$\Rplus\times\Rn$.  Suppose further that $v(0,x)=0$, $\partial_t
v(0,x)=g(x)$, and supp $g(x)\subset \{|x|<1\}$.  Then, for any
$N=0,1,2,...$,
$$\|v'(t,\cdot)\|_{L^2(|x|<t/2)}\le C (1+t)^{-n/2}
\|g\|_{L^2(\R^n)}.$$
\end{lemma}

\noindent{\bf Proof of Lemma \ref{prop2lemma}:} For $t<3$, the
lemma follows easily from \eqref{stenergy}.  We will, thus, assume
$t\ge 3$.  We will show how to get the bound for $\partial_t v$. A
similar argument can be used to bound $\partial_j v$ for
$j=1,2,...,n$.  Since, by Holder's inequality,
$$\|\partial_t v(t,\cdot)\|_{L^2(|x|<t/2)}\le C
t^{n/2}\|\partial_t v(t,\cdot)\|_{L^\infty(|x|<t/2)},$$ it will
suffice to show
\begin{equation}\label{inftyGoal}
\|\partial_t v(t,\cdot)\|_{L^\infty(|x|<t/2)}\le C
t^{-n}\|g\|_{L^2(\R^n)}.
\end{equation}

Since $\partial_t v$ is a linear combination of $e^{\pm
it\sqrt{-\Delta}}g$, it will suffice to show the bound for
$e^{it\sqrt{-\Delta}}g$.  The other piece will follow from the
same argument.

We begin by fixing a smooth, radial cutoff $\chi$ such that
$\chi(\xi)\equiv 1$ for $|\xi|\le 1$ and $\chi(\xi)\equiv 0$ for
$|\xi|\ge 2$.  Then, set
$$\beta(\xi)=\chi(\xi)-\chi(2\xi).$$
Thus, supp $\beta\subset \{1/2 \le |\xi|\le 2\}$, and we have a
partition of unity
$$\chi(\xi)+\sum_{j=1}^\infty \beta(\xi/2^j)=1$$
for all $\xi\ne 0$.  We can then decompose the left side of
\eqref{inftyGoal} as
\begin{multline}\label{decomposed}
\|e^{it\sqrt{-\Delta}}g\|_{L^\infty(|x|<t/2)}\\\le
\|e^{it\sqrt{-\Delta}}\chi(\sqrt{-\Delta})g\|_{L^\infty(|x|<t/2)}
+\sum_{j=0}^\infty\|e^{it\sqrt{-\Delta}}\beta(\sqrt{-\Delta}/2^j)g\|_{L^\infty(|x|<t/2)}
\end{multline}
and examine the pieces on the right side separately.

Since $g$ is supported in $\{|x|\le 1\}$, we see that
\begin{equation}\begin{split}\label{first}
|e^{it\sqrt{-\Delta}}\beta(\sqrt{-\Delta}/2^j)g|&\le \int
\Bigl|\int e^{i(x-y)\cdot\xi}e^{it|\xi|}\beta(\xi/2^j)\:d\xi\Bigr|
|g(y)|\:dy\\
&\le \sup_{|y|\le 1} \left(\left|\int e^{i(x-y)\cdot\xi}
e^{it|\xi|}\beta(\xi/2^j)\:d\xi\right|\right)\:\|g\|_{L^2(\R^n)},
\end{split}
\end{equation}
and similarly
\begin{equation}\label{second}
|e^{it\sqrt{-\Delta}}\chi(\sqrt{-\Delta})g|\le
 \sup_{|y|\le 1} \left(\left|\int e^{i(x-y)\cdot\xi}
e^{it|\xi|}\chi(\xi)\:d\xi\right|\right)\:\|g\|_{L^2(\R^n)}.\end{equation}

For \eqref{first}, if we write the kernel in polar coordinates and
do a change of variables, we have
\begin{align*}
\int e^{i(x-y)\cdot\xi} e^{it|\xi|}\beta(\xi/2^j)\:d\xi&\le
\int_0^\infty
\int_{S^{n-1}}e^{i\rho[(x-y)\cdot\omega+t]}\beta(\rho/2^j)\rho^{n-1}\:d\sigma(\omega)\,d\rho\\
&\le
2^{jn}\int_0^{\infty}\int_{S^{n-1}}e^{i2^j\rho[(x-y)\cdot\omega+t]}a(\rho)\:d\sigma(\omega)\,d\rho
\end{align*}
where $a(\rho)$ is the smooth function, compactly supported away
from 0 given by $\beta(\rho)\rho^{n-1}$.  If we set
$$I_j=2^{jn}\int_0^\infty e^{i\rho[2^j(x-y)\cdot\omega + 2^j
t]}a(\rho)\:d\rho$$ and integrate by parts $N$ times, we see that
$$|I_j|\le C\frac{2^{jn}}{2^{jN}\bigl|t-|x-y|\bigr|^N}\le C
2^{j(n-N)} t^{-N}$$ on $\{|x|<t/2\}\cap\{t\ge 3\}$.  Thus, if we
choose $N>n$, we see that
\begin{equation}\label{firstA}
|e^{it\sqrt{-\Delta}}\beta(\sqrt{-\Delta}/2^j)g|\le C
2^{-jm}t^{-n} \|g\|_{L^2(\R^n)}
\end{equation}
 for some $m>0$.

For \eqref{second}, if we write the kernel in polar coordinates,
we have
$$\int e^{i(x-y)\cdot\xi} e^{it|\xi|}\chi(\xi)\:d\xi\le
\int_0^{\infty}\int_{S^{n-1}}e^{i\rho[(x-y)\cdot\omega+t]}a_0(\rho)
\:d\sigma(\omega)\,d\rho$$ where $a_0$ is the smooth function
given by $\chi(\rho)\rho^{n-1}$.  Here
$\frac{\partial^N}{\partial^N\rho}a_0=0$ for $N<n-1$.  Thus, if we
set
$$I_0=\int_0^\infty
e^{i\rho[(x-y)\cdot\omega+t]}a_0(\rho)\:d\rho,$$ we can integrate
by parts $n$ times to get
$$|I_0|\le
\left|\frac{1}{(t+\omega\cdot(x-y))^n}\right|+\left|\int_0^\infty
\frac{e^{i\rho[(x-y)\cdot\omega+t]}}{(t+\omega\cdot(x-y))^n}a_0^{(n)}(\rho)\:d\rho\right|\le
C t^{-n}$$ on the set $\{|x|<t/2\}\cap\{t\ge 3\}$.  Substituting
this into \eqref{second}, we have
\begin{equation}\label{secondA}
|e^{it\sqrt{-\Delta}}\chi(\sqrt{-\Delta})g|\le C t^{-n}
\|g\|_{L^2(\R^n)}.\end{equation}

Plugging \eqref{firstA} and \eqref{secondA} into
\eqref{decomposed} yields \eqref{inftyGoal} as desired.
 \qed

\noindent{\bf Proof of Proposition \ref{prop2}:} Since the vector
fields $Z$ commute with $\square$ and since
$[Z,\partial_j]=\sum_{i=0}^n a_{ij}\partial_i$ for some constants
$a_{ij}$, it will suffice to show the result for $N=0$.

Let $G_j(s,x)=\chi_{[j,j+1]}(s)G(s,x)$ where $\chi_{[j,j+1]}$ is
the characteristic function of the interval $[j,j+1]$.  Then, let
$v_j$ be the forward solution to $\Box v_j=G_j$ with vanishing
Cauchy data.  By the Cauchy-Schwartz inequality and finite
propagation speed, we have
$$v=\sum_{j=0}^\infty v_j\le C\left(\sum_{j=0}^\infty
|(s-j-|x|+2)v_j|^2\right)^{1/2}.$$  Thus, by Minkowski's integral
inequality,
\begin{equation}\begin{split}\label{breakdown}
\|(1+r)^{-(n-1)/4}&v'(s,\cdot)\|^2_{L^2(\R^n)}\\
&\le C\sum_j
\|(1+r)^{-(n-1)/4}(s-j-|x|+2)v'_j(s,\cdot)\|^2_{L^2(\R^n)}\\
&\le C\sum_j
(s-j+1)^{-(n-1)/2}\|(s-j-|x|+2)v'_j(s,\cdot)\|^2_{L^2(|x|>(s-j)/2)}\\
&\quad\quad\quad\quad+C\sum_j
(s-j+2)^2\|v'_j(s,\cdot)\|^2_{L^2(|x|<(s-j)/2)}.
\end{split}\end{equation}

By Lemma \ref{weightedEnergy} and Duhamel's principle, the first
term in the right side of \eqref{breakdown} is bounded by
\begin{multline*}
(s-j+1)^{-(n-1)/2}\left(\int_j^{j+1}
\|G(\tau,\cdot)\|_{L^2(\R^n)}\:d\tau\right)^2 \\ \le C
\int_j^{j+1}(s-\tau+1)^{-(n-1)/2}\|G(\tau,\cdot)\|^2_{L^2(\R^n)}\:d\tau.
\end{multline*}
By Lemma \ref{prop2lemma} and Duhamel's principle, the second term
in the right side of \eqref{breakdown} is bounded by
\begin{multline*}
(s-j+1)^2\left(\int_j^{j+1} (s-\tau+1)^{-n/2}
\|G(\tau,\cdot)\|_{L^2(\R^n)}\:d\tau\right)^2 \\ \le C
\int_j^{j+1}
(s-\tau+1)^{-n+2}\|G(\tau,\cdot)\|^2_{L^2(\R^n)}\:d\tau.\end{multline*}
Thus, after summing, we have
\begin{multline}\label{breakdown2}
\|(1+r)^{-(n-1)/4}v'(s,\cdot)\|^2_{L^2(\R^n)} \le C \int
(s-\tau+1)^{-(n-1)/2}\|G(\tau,\cdot)\|^2_{L^2(\R^n)}\:d\tau \\+ C
\int (s-\tau+1)^{-n+2}\|G(\tau,\cdot)\|^2_{L^2(\R^n)}\:d\tau.
\end{multline}
If we integrate both sides of \eqref{breakdown2} from $0$ to $t$
and apply Young's inequality, we see that the result follows for
$n\ge 4$. \qed

The next estimate is an $L^1_tL^2_x\to L^2_tL^2_x$ estimate when
the forcing term is not assumed to be compactly supported in $x$.
This follows from the arguments of Smith and Sogge \cite{SS} and
the author \cite{Metcalfe}.  From \cite{Metcalfe}, we have

\begin{lemma}\label{lemma2}
Let $\beta$ be a smooth function supported in $\{|x|<4\}$.  Then,
$$\int_{-\infty}^{\infty}\|\beta(\cdot)e^{it\sqrt{-\Delta}}f(\cdot)\|^2_{\dot{H}^1(\Rn)}\:dt
\le C_{\beta}\|f\|^2_{\dot{H}^1(\Rn)}.$$
\end{lemma}

\noindent From this, we are able to deduce the following corollary
using Duhamel's principle.

\begin{corr}\label{lemma3}
Let $v$ be a solution to the wave equation $\square v=G$ in
$\Rplus \times \Rn$.  Then, for $\beta$ a smooth function
supported in $\{|x|\le 4\}$, we have
$$\sup_{|\alpha|\le 1} \|\beta(x) \partial^\alpha_{t,x}v(s,x)\|_{L^2_{s,x}([0,t]\times\Rn)}\le
C\|v'(0,\cdot)\|_{L^2(\Rn)}+C\int_0^t
\|G(s,\cdot)\|_{L^2(\Rn)}\:ds.$$
\end{corr}

And, from Corollary \ref{lemma3}, we then easily get the following
corollary.
\begin{corr}\label{cor5}
Let $v$ be a solution to the wave equation $\square v=G$ in
$\Rplus \times \Rn$.  Then, for $\beta$ a smooth function
supported in $\{R\le |x|\le 4R\}$ with $R>1$, we have
$$\|r^{-1/2} \beta(x) v'(s,x)\|_{L^2_{s,x}([0,t]\times\Rn)}\le
C\|v'(0,\cdot)\|_{L^2(\Rn)}+C\int_0^t
\|G(s,\cdot)\|_{L^2(\Rn)}\:ds.$$
\end{corr}

\noindent{\bf Proof of Corollary \ref{cor5}:} As above, setting
$v_R(t,x)=v(Rt,Rx)$, $G_R(t,x)=R^2G(Rt,Rx)$, and
$\beta_R(x)=\beta(Rx)$, we have $\square v_R=G_R$ and $\supp
\beta_R\subset \{1\le |x|\le 4\}$.  Thus, by Lemma \ref{lemma3},
\begin{align*}
\|r^{-1/2}\beta_R(x)v_R'(s,x)\|_{L^2_{s,x}([0,t/R]\times \Rn)}&\le
\|\beta_R(x)v_R'(s,x)\|_{L^2_{s,x}([0,t/R]\times\Rn)}\\
&\le C\|v_R'(0,\cdot)\|_{L^2(\Rn)}+C\int_0^{t/R}
\|G_R(s,\cdot)\|_{L^2_x(\Rn)}\:ds
\end{align*}
or \begin{multline*}
R\|r^{-1/2}\beta(Rx)v'(Rs,Rx)\|_{L^2_{s,x}([0,t/R]\times\Rn)}\\\le
CR\|v'(0,Rx)\|_{L^2_x(\Rn)}+CR^2\int_0^{t/R}\|G(Rs,Rx)\|_{L^2_x(\Rn)}\:ds.\end{multline*}
After a change of variables, this becomes
$$\|r^{-1/2} \beta(x) v'(s,x)\|_{L_{s,x}([0,t]\times\Rn)}\le
C\|v'(0,\cdot)\|_{L^2(\Rn)}+C\int_0^t
\|G(s,\cdot)\|_{L^2(\Rn)}\:ds$$ as desired.\qed

We are now ready to prove our key radial decay estimate.

\begin{proposition}\label{prop6}  Suppose $n\ge 4$.  Let
$v$ be a solution of the wave equation $\square v=G$ in
$\Rplus\times\Rn$.  Then,
$$\|(1+r)^{-(n-1)/4}v'(s,x)\|_{L^2_{s,x}([0,t]\times\Rn)}
\le
C\|v'(0,\cdot)\|_{L^2(\Rn)}+C\int_0^t\|G(s,\cdot)\|_{L^2_x(\Rn)}\:ds
$$
for a uniform constant $C$.
\end{proposition}

\noindent{\bf Proof of Proposition \ref{prop6}:} Fix $\chi$ to be
a smooth, radial cutoff function where $\chi(x)\equiv 1$ when
$|x|\le 1$ and $\chi(x)\equiv 0$ for $|x|>2$.  Additionally, set
$\beta(x)=\chi(x)-\chi(2x)$.  Thus, $\supp\beta\subset\{1/2\le
|x|\le 2\}$ and
\begin{align*}
\|(1+r)^{-(n-1)/4}&v'(s,x)\|^2_{L^2_{s,x}([0,t]\times
\Rn)}\\
&\le \|(1+r)^{-(n-1)/4}\chi(x)v'(s,x)\|^2_{L^2_{s,x}([0,t]\times
\Rn)} \\&\quad\quad\quad+ \sum_{j=0}^{\infty}
\|(1+r)^{-(n-1)/4}\beta(x/2^j)v'(s,x)\|^2_{L^2_{s,x}([0,t]\times
\Rn)}\\
&\le \|\chi(x)v'(s,x)\|^2_{L^2_{s,x}([0,t]\times \Rn)}\\
&\quad\quad\quad +\sum_{j=0}^{\infty}
2^{-j(n-3)/4}\|(1+r)^{-1/2}\beta(x/2^j)v'(s,x)\|^2_{L^2_{s,x}([0,t]\times
\Rn)}.
\end{align*}
The result, then, follows from an application of Corollary
\ref{lemma3} and Corollary \ref{cor5}.\qed

Since the vector fields $Z$ commute with $\square$ and since
$[Z,\partial_j]=\sum_{i=0}^n a_{ij}\partial_i$ for some constants
$a_{ij}$, Proposition \ref{prop6} and \eqref{stenergy} imply

\begin{theorem}\label{thm7}
Suppose that $v$ is a solution of the wave equation $\square v=G$
in $\Rplus\times\Rn$.  Then, for any $N=0,1,2,...$
\begin{multline*}
\sum_{|\alpha|\le N} \left(\|Z^\alpha v'(t,\cdot)\|_{L^2_x(\Rn)} +
\|(1+r)^{-(n-1)/4}Z^\alpha v'(s,x)\|_{L^2_{s,x}([0,t]\times\Rn)}\right)\\
\le C\sum_{|\alpha|\le N}\|Z^\alpha
v'(0,\cdot)\|_{L^2(\Rn)}+C\sum_{|\alpha|\le N}\int_0^t\|Z^\alpha
G(s,\cdot)\|_{L^2_x(\Rn)}\:ds
\end{multline*}
\end{theorem}

In addition to the $L^2$ estimates, we will need the following
pointwise estimate.  This is a weighted Sobolev estimate (see
e.g., \cite{KSS2}).

\begin{lemma}\label{lemmaii.1}
Suppose that $h\in C^\infty(\Rn)$.  Then, for $R>1$,
$$\|h\|_{L^\infty(R/2\le |x|\le R)}\le CR^{-(n-1)/2}\sum_{|\alpha|\le (n+2)/2}
\|Z^\alpha h\|_{L^2(R/4\le |x|\le 2R)}.$$
\end{lemma}

\noindent{\bf Proof of Lemma \ref{lemmaii.1}:} By Sobolev's lemma
for $\mathbb{R}\times S^{n-1}$, we have
$$|h(x)|\le C \sum_{|\alpha|+j\le
\frac{n+2}{2}}\left(\int_{|x|-1/4}^{|x|+1/4}\int_{S^{n-1}}|\partial^j_r\Omega^\alpha
h(r\omega)|^2\:d\sigma(\omega)\:dr\right)^{1/2}.$$ Thus,
\begin{align*}
\|h(x)\|_{L^\infty(R/2\le |x|\le R)} &\le C R^{-(n-1)/2}
\sum_{|\alpha|+j\le
\frac{n+2}{2}}\left(\int_{R/4}^{2R}\int_{S^{n-1}}|\partial^j_r\Omega^\alpha
h(r\omega)|^2 r^{n-1} \:d\sigma(\omega)\:dr\right)^{1/2}\\
&\le C R^{-(n-1)/2}\sum_{|\alpha|\le \frac{n+2}{2}} \|Z^\alpha
h\|_{L^2(R/4\le |x|\le 2R)}.
\end{align*}\qed


\newsection{Global existence in Minkowski space.}
We now want to use Theorem \ref{thm7} and Lemma \ref{lemmaii.1} to
prove Theorem \ref{GEminkowski}.  We will use an iteration to
solve \eqref{mainMinkowski} and to show that
\begin{equation}\label{bounded}
\sup_{0\le t\le T} \sum_{|\alpha|\le n+2} \Bigl(\|Z^\alpha
u'(t,\cdot)\|_{L^2(\Rn)} +\|(1+r)^{-(n-1)/4}Z^\alpha
u'\|_{L^2_{s,x}([0,T]\times\Rn)}\Bigr)\le C\varepsilon
\end{equation}
for any time $T$.

\noindent{\bf Proof of Theorem \ref{GEminkowski}:} Set
$u_{-1}\equiv 0$.  Then, define $u_k$ recursively by setting it to
be the solution of
\begin{equation}\label{recursion}
\begin{cases}
\Box u_k(t,x)=Q(u'_{k-1}(t,x)), \quad (t,x)\in [0,T_*]\times \Rn\\
u_k(0,\cdot)=f,\quad \partial_t u_k(0,\cdot)=g.
\end{cases}
\end{equation}
Let
$$M_k(T)=\sup_{0\le t\le T} \sum_{|\alpha|\le n+2}
\|Z^\alpha u'_k(t,\cdot)\|_{L^2(\Rn)} +\sum_{|\alpha|\le
n+2}\|(1+r)^{-(n-1)/4}Z^\alpha u'_k\|_{L^2_{s,x}([0,T]\times\Rn)}.
$$

By \eqref{GEminkowskiAssumption} and Theorem \ref{thm7}, there is
a constant $C_0$ such that
$$M_0(T)\le C_0\varepsilon$$
for all $T$.  Our goal is to inductively prove that if
$\varepsilon < \varepsilon_0$ is sufficiently small, then
\begin{equation}\label{InductiveGoal}
M_k(T)\le 2C_0\varepsilon
\end{equation}
 for every $k=1,2,3...$.  To do so, we
assume that this bound holds for $k-1$ and we will use the
assumption to prove \eqref{InductiveGoal} for $k$.  By Theorem
\ref{thm7}, we have
$$M_k(T)\le C_0\varepsilon + C\sum_{|\alpha|\le n+2} \int_0^T
\|Z^\alpha Q(u'_{k-1})(s,\cdot)\|_{L^2(\Rn)}\:ds.$$ Since $Q$ is
quadratic, the bound
$$|Z^\alpha Q(u'_{k-1})(s,x)|\le C\left(\sum_{|\alpha|\le n+2}
|Z^\alpha u'_{k-1}(s,x)|\right)\left(\sum_{|\alpha|\le
(n+2)/2}|Z^\alpha u'_{k-1}(s,x)|\right)$$ holds for all
$|\alpha|\le n+2$.  Thus, applying Lemma \ref{lemmaii.1}, we have
\begin{align*}
&\sum_{|\alpha|\le n+2} \|Z^\alpha
Q(u'_{k-1})(s,\cdot)\|_{L^2(\{|x|\in [2^j,2^{j+1}]\})}\\
&\le C \sum_{|\alpha|\le n+2}\|Z^\alpha
u'_{k-1}(s,x)\|_{L^2(\{|x|\in [2^j,2^{j+1}]\})}\sum_{|\alpha|\le
(n+2)/2}\|Z^\alpha
u'_{k-1}(s,x)\|_{L^\infty(\{|x|\in [2^j,2^{j+1}]\})}\\
&\le C2^{-j(n-1)/2}\sum_{|\alpha|\le n+2}\|Z^\alpha
u'_{k-1}(s,x)\|^2_{L^2(\{|x|\in [2^{j-1},2^{j+2}]\})}\\
&\le C \sum_{|\alpha|\le n+2}\|(1+r)^{-(n-1)/4}Z^\alpha
u'_{k-1}(s,x)\|^2_{L^2(\{|x|\in [2^{j-1},2^{j+2}]\})}
\end{align*}
By the standard Sobolev lemma, we also have
$$\sum_{|\alpha|\le n+2} \|Z^\alpha
Q(u'_{k-1})(s,\cdot)\|_{L^2(\{|x|< 1\})}\le \sum_{|\alpha|\le
n+2}\|Z^\alpha u'_{k-1}(s,\cdot)\|^2_{L^2(\{|x|< 2\})}.
$$
Hence,
$$\sum_{|\alpha|\le n+2} \|Z^\alpha Q(u'_{k-1})(s,\cdot)\|_{L^2(\Rn)}\le C
\sum_{|\alpha|\le n+2} \|(1+r)^{-(n-1)/4}Z^\alpha
u'_{k-1}(s,\cdot)\|^2_{L^2(\Rn)}$$ and using the inductive
hypothesis,
\begin{align*}
M_k(T)&\le C_0\varepsilon + C\sum_{|\alpha|\le n+2} \int_0^T
\|Z^\alpha Q(u'_{k-1})(s,\cdot)\|_{L^2(\Rn)}\:ds\\
&\le C_0\varepsilon + C\sum_{|\alpha|\le n+2}
\|(1+r)^{-(n-1)/4}Z^\alpha
u'_{k-1}(s,\cdot)\|^2_{L^2(\{(s,x): 0\le s\le T\})}\\
&\le C_0\varepsilon + CM_{k-1}^2(T)\\
&\le C_0\varepsilon + 4C_0^2C\varepsilon^2.
\end{align*}
Thus, if $\varepsilon\le \varepsilon_0=\frac{1}{8CC_0}$, then we
see that \eqref{InductiveGoal} holds for any $T$.
\par
We now need to show that the $u_k$ converge to a solution. If we
set
\begin{multline}\label{ak}
A_k(T)=\sup_{0\le t\le T} \sum_{|\alpha|\le n+2} \|Z^\alpha
(u'_k-u'_{k-1})(t,\cdot)\|_{L^2(\Rn)}\\
+\sum_{|\alpha|\le n+2}\|(1+r)^{-(n-1)/4}Z^\alpha
(u'_k-u'_{k-1})\|_{L^2(\{(s,x):0\le s\le T\})},
\end{multline}
the proof will be completed if we can show
\begin{equation}\label{ConvergenceGoal}
A_k(T)\le \frac{1}{2}A_{k-1}(T), \quad k=1,2,3,...
\end{equation}
for any $T$. Since $Q$ is quadratic, we have
$$|Q(u'_{k-1})-Q(u'_{k-2})|\le
C\left(|u'_{k-1}||u'_{k-1}-u'_{k-2}|+|u'_{k-2}||u'_{k-1}-u'_{k-2}|\right).$$
By repeating the previous arguments, we have
\begin{multline*}
A_k(T)\le C\sum_{|\alpha|\le n+2}\int_0^T
\|(1+r)^{-(n-1)/4}Z^\alpha
(u'_{k-1}-u'_{k-2})(s,\cdot)\|_{L^2(\Rn)} \\
\times \left(\|(1+r)^{-(n-1)/4}Z^\alpha
u'_{k-1}(s,\cdot)\|_{L^2(\Rn)}+\|(1+r)^{-(n-1)/4}Z^\alpha
u'_{k-2}(s,\cdot)\|_{L^2(\Rn)}\right)\: ds.
\end{multline*}
Thus, by the Schwarz inequality, we have
$$A_k(T)\le C[M_{k-1}(T)+M_{k-2}(T)]A_{k-1}(T)\le 4CC_0\varepsilon A_{k-1}(T)$$
for any $T$.  For $\varepsilon$ as above, we see that
\eqref{ConvergenceGoal} holds which completes the proof. \qed


\newsection{Main estimates in the exterior domain.}

In the next section, we will prove Theorem \ref{GEexterior}. The
first step in adapting the argument of the preceding section is to
show that there are analogs of Theorem \ref{thm7} and Lemma
\ref{lemmaii.1} that hold exterior to a nontrapping obstacle.
\par
By scaling, we may assume that $\mathcal{K}\subset\{|x|<1/2\}$.
\par
The analog of Lemma \ref{lemmaii.1} follows directly from the
proof given above when $R>2$. That is, if $h(x)=0$ when $x\in
\partial\mathcal{K}$, then
\begin{equation}\label{inftyEstimate}
\|h\|_{L^\infty(R/2\le |x|\le R)}\le C
R^{-(n-1)/2}\sum_{|\alpha|\le (n+2)/2} \|Z^\alpha h\|_{L^2(R/4\le
|x|\le 2R)}.
\end{equation}
When $R\le 2$, this follows from standard Sobolev estimates.
\par
We will also need exterior domain analogs of the energy-type
estimates.  In order to avoid issues with compatibility
conditions, we will restrict to the case where the initial data
vanish.  In proving Theorem \ref{GEexterior}, we will reduce to
this situation.  Thus, we will be looking at solutions of
\begin{equation}\label{ExteriorEquation}
\begin{cases}
\Box w(t,x)=F(t,x),\quad (t,x)\in \Rplus\times \Omega\\
w(0,x)=\partial_t w(0,x)=0\\
w(t,x)=0 \quad\text{for } x\in \partial\Omega.
\end{cases}
\end{equation}
\par
It is well known that analog of \eqref{stenergy} holds in domains
exterior to a nontrapping obstacle.  Specifically, if $w$ is given
by \eqref{ExteriorEquation}, then
\begin{equation}\label{stenergyExterior}
\|w'(t,\cdot)\|_{L^2(\Omega)}\le C\int_0^t
\|F(s,\cdot)\|_{L^2(\Omega)}\:ds.
\end{equation}

We now turn our attention to the exterior domain analog of our
main estimate.

\begin{theorem}\label{mainExterior}
Suppose $n\ge 4$, and let $w$ be the solution of
\eqref{ExteriorEquation}. Then, for any $N=0,1,2,...$
\begin{multline*}
\sum_{|\alpha|\le N}\left(\|Z^\alpha
w'(t,\cdot)\|_{L^2(\Omega)}+\|(1+r)^{-(n-1)/4}Z^\alpha
w'\|_{L^2_{s,x}([0,t]\times\Omega)}\right)\\
\le C\sum_{|\alpha|\le N}\int_0^t \|Z^\alpha
F(s,\cdot)\|_{L^2(\Omega)}\:ds + C\sup_{0\le s\le
t}\sum_{|\alpha|\le N-1}\|Z^\alpha F(s,\cdot)\|_{L^2(\Omega)}\\
+C\sum_{|\alpha|\le N-1}\|Z^\alpha
F\|_{L^2_{s,x}([0,t]\times\Omega)}.
\end{multline*}
\end{theorem}

For any $R>1/2$, set $B_R=\{|x|\le R\}\cap \Omega$ and
$E_R=\{|x|\ge R\}$. We will prove Theorem \ref{mainExterior} by
proving the following four estimates:
\begin{multline}\label{mainExterior1}
\sum_{|\alpha|\le N}\|Z^\alpha w'(t,\cdot)\|_{L^2(B_2)} \le
C\sum_{|\alpha|\le N}\int_0^t \|Z^\alpha
F(s,\cdot)\|_{L^2(\Omega)}\:ds \\+ C\sup_{0\le s\le
t}\sum_{|\alpha|\le N-1}\|Z^\alpha F(s,\cdot)\|_{L^2(\Omega)},
\end{multline}

\begin{multline}\label{mainExterior2}
\sum_{|\alpha|\le N}\|(1+r)^{-(n-1)/4}Z^\alpha
w'\|_{L^2_{t,x}([0,t]\times B_2)} \le C\sum_{|\alpha|\le
N}\int_0^t \|Z^\alpha F(s,\cdot)\|_{L^2(\Omega)}\:ds \\+
C\sum_{|\alpha|\le N-1}\|Z^\alpha
F\|_{L^2_{s,x}([0,t]\times\Omega)},
\end{multline}

\begin{multline}\label{mainExterior3}
\sum_{|\alpha|\le N}\|Z^\alpha w'(t,\cdot)\|_{L^2(E_2)} \le
C\sum_{|\alpha|\le N}\int_0^t \|Z^\alpha
F(s,\cdot)\|_{L^2(\Omega)}\:ds \\
+C\sum_{|\alpha|\le N-1}\|Z^\alpha
F\|_{L^2_{s,x}([0,t]\times\Omega)},
\end{multline}

\begin{multline}\label{mainExterior4}
\sum_{|\alpha|\le N}\|(1+r)^{-(n-1)/4}Z^\alpha
w'\|_{L^2_{t,x}([0,t]\times E_2)} \le C\sum_{|\alpha|\le
N}\int_0^t \|Z^\alpha F(s,\cdot)\|_{L^2(\Omega)}\:ds\\ +
C\sum_{|\alpha|\le N-1}\|Z^\alpha
F\|_{L^2_{s,x}([0,t]\times\Omega)}.
\end{multline}

\noindent{\bf Proof of Equation \eqref{mainExterior1}:} Since
$$\sum_{|\alpha|\le N}\|Z^\alpha w'(t,\cdot)\|_{L^2(B_2)}\le
C\sum_{|\alpha|\le
N}\|\partial^{\alpha}_{t,x}w'(t,\cdot)\|_{L^2(\Omega)},$$ we need
only show
\begin{equation}\label{mainExterior1a}
\sum_{|\alpha|\le
N}\|\partial^\alpha_{t,x}w'(t,\cdot)\|_{L^2(\Omega)}\le C\int_0^t
\sum_{|\alpha|\le N}\|\partial^{\alpha}_{t,x}
F(s,\cdot)\|_{L^2(\Omega)}+C\sum_{|\alpha|\le
N-1}\|\partial^\alpha_{t,x}F(t,\cdot)\|_{L^2(\Omega)}.
\end{equation}
We will prove \eqref{mainExterior1a} via induction.  When $N=0$,
this follows from the standard energy inequality
\eqref{stenergyExterior}.  Thus, we will assume that
\eqref{mainExterior1a} holds for $N-1$ and prove that this implies
the result for $N$.
\par Notice that since $\partial_t$ preserves the boundary
condition and that $[\Box,\partial_t]=0$, the inductive hypothesis
applied to $\partial_t w$ gives
\begin{multline}\label{indHyp}
\sum_{|\alpha|\le N-1}\|\partial^\alpha_{t,x}\partial_t
w'(t,\cdot)\|_{L^2(\Omega)}\\
\le C\sum_{|\alpha|\le N-1}\int_0^t \|\partial_{t,x}^\alpha
\partial_t F(s,\cdot)\|_{L^2(\Omega)}\:ds +C\sum_{|\alpha|\le
N-2}\|\partial_{t,x}^\alpha \partial_t F(t,\cdot)\|_{L^2(\Omega)}.
\end{multline}
Thus, since $\|\partial^\alpha_{t,x}\partial_t^2
w(t,\cdot)\|_{L^2(\Omega)}\le \|\partial^\alpha_{t,x}\partial_t
w'(t,\cdot)\|_{L^2(\Omega)}$ and $\partial_t^2w=\Delta w+F$, we
have
\begin{multline}\label{indHypLaplace}
\sum_{|\alpha|\le N-1}\|\partial_{t,x}^\alpha \Delta
w(t,\cdot)\|_{L^2(\Omega)}\\
\le C\sum_{|\alpha|\le N}\int_0^t \|\partial^\alpha_{t,x}
F(s,\cdot)\|_{L^2(\Omega)}\:ds + C\sum_{|\alpha|\le
N-1}\|\partial^\alpha_{t,x} F(t,\cdot)\|_{L^2(\Omega)}.
\end{multline}

Using elliptic regularity, we see that
\begin{align*}
&\sum_{|\alpha|\le
N}\|\partial^\alpha_{t,x}w'(t,\cdot)\|_{L^2(\Omega)}\\
&\le \|w'(t,\cdot)\|_{L^2(\Omega)} +\sum_{1\le |\alpha|\le
N}\|\partial_x^\alpha w'(t,x)\|_{L^2(\Omega)}+\sum_{|\alpha|\le
N-1}\|\partial^{\alpha}_{t,x}\partial_t w'(t,x)\|_{L^2(\Omega)}\\
&\le C\sum_{|\alpha|\le N-1}\left(\|\partial_x^\alpha
w'(t,\cdot)\|_{L^2(\Omega)}+\|\partial_x^\alpha \Delta
w(t,x)\|_{L^2(\Omega)}+\|\partial^{\alpha}_{t,x}\partial_t
w'(t,x)\|_{L^2(\Omega)}\right)
\end{align*}

Thus, by the inductive hypothesis, \eqref{indHyp} and
\eqref{indHypLaplace}, we see that the proof of
\eqref{mainExterior1} is complete. \qed

Since
$$
\sum_{|\alpha|\le N}\|(1+r)^{-(n-1)/4}Z^\alpha
w'\|_{L^2_{s,x}([0,t]\times B_2)}\le C\sum_{|\alpha|\le
N}\|\partial^\alpha_{t,x} w'\|_{L^2_{s,x}([0,t]\times B_2)},$$ we
see that \eqref{mainExterior2} follows easily from the following
lemma.

\begin{lemma}\label{mainExterior2lemma}
Let $w$ be as in \eqref{ExteriorEquation}.  Then, for any
$N=0,1,2,...$,
\begin{multline*}
\sum_{|\alpha|\le N}\|\partial_{t,x}^\alpha
w'\|_{L^2_{s,x}([0,t]\times B_2)}\\
\le C\sum_{|\alpha|\le N}\int_0^t \|\partial_{t,x}^\alpha
F(s,\cdot)\|_{L^2(\Omega)}\:ds + C\sum_{|\alpha|\le
N-1}\|\partial_{t,x}^\alpha F\|_{L^2_{s,x}([0,t]\times\Omega)}
\end{multline*}
\end{lemma}

\noindent{\bf Proof of Lemma \ref{mainExterior2lemma}:} Following
the same induction argument as above, it will suffice to show the
$N=0$ case,
$$\|w'\|_{L^2_{s,x}([0,t]\times B_2)}\le C\int_0^t
\|F(s,\cdot)\|_{L^2(\Omega)}\:ds.$$
\par
Suppose that $F(s,x)=0$ when $|x|>4$.  In this case, since we are
exterior to a nontrapping obstacle, the local energy decay of
Melrose \cite{Melrose} and Duhamel's principle imply
$$\|w'(t,\cdot)\|_{L^2(B_2)}\le C \int_0^t [1+(t-s)]^{-n/2}
\|F(s,\cdot)\|_{L^2(\Omega)}\:ds$$ Thus, from Minkowski's integral
inequality, we have
$$\|w'\|_{L^2_{s,x}([0,t]\times B_2)}\le C\int_0^t
\|F(s,\cdot)\|_{L^2(\Omega)}\:ds$$ as desired.
\par
Now suppose that $F(s,x)=0$ for $|x|\le 4$.  Fix $\rho\in
C_0^\infty$ such that $\rho\equiv 1$ when $|x|\le 2$ and
$\rho\equiv 0$ when $|x|>4$.  Let $u_0$ be the solution to the
free wave equation $\Box u_0=F$ with vanishing data.  Here we have
set $F$ to zero on $\Rn \backslash \Omega$.  Write
$$w=u_0+u_r=(1-\rho)u_0 + [\rho u_0+u_r]$$
and set $v=\rho u_0+u_r$.  By our assumption on $F$, we have
$$\Box v = \nabla_x \rho \cdot \nabla u_0 + (\Delta \rho) u_0$$
The argument in the preceding paragraph implies
$$\|w'\|_{L^2_{s,x}([0,t]\times
B_2)}=\|v'\|_{L^2_{s,x}([0,t]\times B_2)}\le \|\rho'
u'_0\|_{L^2_{s,x}([0,t]\times \Rn)}+\|(\Delta \rho)
u_0\|_{L^2_{s,x}([0,t]\times \Rn)}.$$ Thus, an application of
Corollary \ref{lemma3} concludes the proof.\qed

\begin{corr}\label{mainExterior2corr}
Let $w$ be as in \eqref{ExteriorEquation}.  Then, for any
$N=0,1,2,...$,
\begin{multline*}
\sum_{|\alpha|\le N}\|\partial_{t,x}^\alpha
w\|_{L^2_{s,x}([0,t]\times B_2)}\\
\le C\sum_{|\alpha|\le N}\int_0^t \|\partial_{t,x}^\alpha
F(s,\cdot)\|_{L^2(\Omega)}\:ds + C\sum_{|\alpha|\le
N-1}\|\partial_{t,x}^\alpha F\|_{L^2_{s,x}([0,t]\times\Omega)}
\end{multline*}
\end{corr}

\noindent{\bf Proof of Corollary \ref{mainExterior2corr}:} By the
previous lemma, it will suffice to show
$$
\|w\|_{L^2_{s,x}([0,t]\times B_2)} \le C\int_0^t \|
F(s,\cdot)\|_{L^2(\Omega)}\:ds
$$
When $F(s,x)=0$ for $|x|>4$, using a modification of local energy
decay, see \cite{Metcalfe}, we have
\begin{align*}
\|w(t,\cdot)\|_{L^2(B_2)}&\le C \int_0^t
(1+t-s)^{-n/2}\|F(s,\cdot)\|_{\dot{H}^{-1}_D(\tilde{\Omega})}\\
&\le C\int_0^t (1+t-s)^{-n/2}\|F(s,\cdot)\|_{L^2(\Omega)}
\end{align*}
 where $\tilde{\Omega}$ is a compact manifold with boundary
containing $B_2$.  Thus, in this case, the result follows from
Young's inequality.

When $F(s,x)=0$ for $|x|\le 4$, we can argue as in the previous
lemma in order to complete the proof of the corollary. \qed

We can now conclude the proof of \eqref{mainExterior} by proving
\eqref{mainExterior3} and \eqref{mainExterior4}.  Let's begin by
fixing a $\beta\in C^\infty$ such that $\beta(x)\equiv 1$ for
$|x|\ge 2$ and $\beta(x)\equiv 0$ for $|x|\le 1$.  Setting
$v=\beta w$, we see that $v=w$ on $|x|>2$ and that $v$ solves the
free wave equation
$$\Box v = \beta F + \nabla\beta\cdot\nabla_x w + (\Delta\beta)
w$$ Decompose $v$ into $v=v_1+v_2$ where $v_1$ solves $\Box
v_1=\beta F$ and $v_2$ is the solution of $\Box v_2=\nabla
\beta\cdot\nabla_x w+(\Delta \beta)w$.  Set $G=\nabla\beta\cdot
\nabla_x w+(\Delta\beta)w$.
\par
By Theorem \ref{thm7}, we have
\begin{multline*}
\sum_{|\alpha|\le N}\left(\|Z^\alpha
v_1'(t,\cdot)\|_{L^2(E_2)}+\|(1+r)^{-(n-1)/4}Z^\alpha
v_1'\|_{L^2_{s,x}([0,t]\times E_2)}\right)\\
\le \sum_{|\alpha|\le N}\int_0^t
\|Z^{\alpha}F(s,\cdot)\|_{L^2(\Omega)}\:ds
\end{multline*}
since for $|\alpha|\le N$,
$$\|Z^{\alpha}\beta(\cdot)F(s,\cdot)\|_{L^2(\Omega)}\le C\sum_{|\alpha|\le
N}\|Z^{\alpha}F(s,\cdot)\|_{L^2(\Omega)}.$$
\par
Thus, it remains to show

\begin{multline}\label{mainExterior3a}
\sum_{|\alpha|\le N}\|Z^\alpha v_2'(t,\cdot)\|_{L^2(E_2)} \le
C\sum_{|\alpha|\le N}\int_0^t \|Z^\alpha
F(s,\cdot)\|_{L^2(\Omega)}\:ds \\
+C\sum_{|\alpha|\le N-1}\|Z^\alpha
F\|_{L^2_{s,x}([0,t]\times\Omega)},
\end{multline}

\begin{multline}\label{mainExterior4a}
\sum_{|\alpha|\le N}\|(1+r)^{-(n-1)/4}Z^\alpha
v_2'\|_{L^2_{t,x}([0,t]\times E_2)} \le C\sum_{|\alpha|\le
N}\int_0^t \|Z^\alpha F(s,\cdot)\|_{L^2(\Omega)}\:ds\\ +
C\sum_{|\alpha|\le N-1}\|Z^\alpha
F\|_{L^2_{s,x}([0,t]\times\Omega)}.
\end{multline}

\noindent{\bf Proof of Equation \eqref{mainExterior3a}:} Let
$G_j(s,x)=\chi_{[j,j+1]}(s)G(s,x)$ where $\chi_{[j,j+1]}$ is the
characteristic function of the interval $[j,j+1]$.  Then, let
$v_{2,j}$ be the forward solution of $\Box v_{2,j}=G_j$ in free
space with zero initial data.  By finite propagation speed and the
Cauchy-Schwartz inequality, we have
\begin{equation}\label{csWeighted}
v_2=\sum_{j=0}^\infty v_{2,j}\le C\left(\sum_{j=0}^\infty
|(t-j-|x|+2)v_{2,j}|^2\right)^{1/2}.
\end{equation}
 Thus, by the Minkowski
integral inequality and Lemma \ref{weightedEnergy},
\begin{align*}
\sum_{|\alpha|\le N}\|Z^{\alpha}v_2'(t,\cdot)\|^2_{L^2(E_2)}&\le
C\sum_{|\alpha|\le N}\sum_{j}\|(t-j-|x|+2)Z^\alpha
v_{2,j}'(t,\cdot)\|^2_{L^2(E_2)}\\
&\le C\sum_{|\alpha|\le N}\sum_{j}\left(\int_j^{j+1}\|Z^\alpha
G(s,\cdot)\|_{L^2(\Rn)}\:ds\right)^2\\
&\le C\sum_{|\alpha|\le N} \|Z^\alpha
G\|^2_{L^2_{s,x}([0,t]\times\Rn)}\\
&\le C\sum_{|\alpha|\le
N}\|\partial^\alpha_{t,x}w'\|^2_{L^2_{s,x}([0,t]\times
B_2)}+C\sum_{|\alpha|\le
N}\|\partial^\alpha_{t,x}w\|^2_{L^2_{s,x}([0,t]\times B_2)}
\end{align*}
Thus, \eqref{mainExterior3a} follows from Lemma
\ref{mainExterior2lemma} and Corollary
\ref{mainExterior2corr}.\qed

\noindent{\bf Proof of Equation \eqref{mainExterior4a}:} By
Proposition \ref{prop2}, we have
\begin{align*}
\sum_{|\alpha|\le N} \|(1+r)^{-(n-1)/4}Z^\alpha
&v'_2\|_{L^2{s,x}([0,t]\times E_2)} \le C \sum_{|\alpha|\le N}
\|Z^\alpha G\|_{L^2_{s,x}(\R^n)}\\
&\le C \sum_{|\alpha|\le
N}\|\partial^\alpha_{t,x}w'\|_{L^2_{s,x}([0,t]\times B_2)}+
C\sum_{|\alpha|\le N}
\|\partial^\alpha_{t,x}w\|_{L^2_{s,x}([0,t]\times B_2)}
\end{align*}
Thus, \eqref{mainExterior4a} follows from Lemma
\ref{mainExterior2lemma} and Corollary
\ref{mainExterior2corr}.\qed


\newsection{Global existence exterior to a nontrapping obstacle.}
We now turn to the proof of Theorem \ref{GEexterior}.  By scaling,
we may assume that the obstacle $\mathcal{K}$ is contained in
$\{|x|<1/2\}$.  It is convenient to show that one can instead show
global existence for an equivalent nonlinear equation which has
vanishing Cauchy data, as in Keel-Smith-Sogge \cite{KSS1}.  This
allows one to avoid the issues regarding the compatibility
conditions. At this point, we can follow an iteration argument
similar to that used to prove Theorem \ref{GEminkowski}.

\noindent{\bf Proof of Theorem \eqref{GEexterior}:} We start by
making the reduction mentioned above.  Notice that if $f,g$
satisfy \eqref{GEexteriorAssumption}, then we can find a local
solution $u$ to \eqref{main} in $0<t<1$.  Moreover, if
$\varepsilon>0$ in \eqref{GEexteriorAssumption} is sufficiently
small, there is a constant $C$ so that
\begin{equation}\label{reduction}
\sup_{0\le t\le 1}\sum_{|\alpha|\le n+2}\|Z^\alpha
u'(t,\cdot)\|_{L^2(\Omega)}+\sum_{|\alpha|\le
n+2}\|(1+r)^{-(n-1)/4}Z^\alpha
u'\|_{L^2_{s,x}([0,1]\times\Omega)}\le C\varepsilon
\end{equation}
To see this, notice that local existence theory (see e.g.,
\cite{KSS1}) implies that \eqref{reduction} holds when
$\varepsilon$ is sufficiently small and the norms on the left side
are taken over $\{|x|<10\}$.  By finite propagation speed, on
$\{0<t<1\}\times \{|x|\ge 10\}$, $u$ agrees with a solution of the
boundaryless wave equation $\Box u=Q(u')$ with data equal to a
cutoff times the original data $(f,g)$.  Thus, in this case,
\eqref{reduction} follows from \eqref{bounded}.
\par
We are now ready to set up the iteration.  Fix $\eta\in
C^\infty(\R)$ such that $\eta(t)\equiv 1$ when $t\le 1/2$ and
$\eta(t)\equiv 0$ for $t>1$.  Let
$$u_0=\eta u.$$
Thus,
$$\Box u_0=\eta Q(u')+[\Box,\eta]u.$$
Hence, in order to show that there is a solution to $\Box u=Q(u')$
for all $t$, it will suffice to show that there is a solution
$w=u-u_0$ of
\begin{equation}\label{wEquation}
\begin{cases}
\Box w=(1-\eta)Q((u_0+w)')-[\Box,\eta](u_0+w)\\
w(t,x)=0\quad\text{for } x\in\partial\Omega\\
w(0,x)=\partial_t w(0,x)=0.
\end{cases}
\end{equation}
\par
In order to set up the iteration, as in the proof of Theorem
\ref{GEminkowski}, set $w_0=0$ and define $w_k$ recursively by
letting it be a solution of
\begin{equation}\label{wRecursiveEquation}
\begin{cases}
\Box w_k=(1-\eta)Q((u_0+w_{k-1})')-[\Box,\eta](u_0+w_k)\\
w_k(t,x)=0\quad\text{for } x\in\partial\Omega\\
w_k(0,x)=\partial_t w_k(0,x)=0.
\end{cases}
\end{equation}
Also, as before, set
$$M_k(T)=\sup_{0\le t\le T} \sum_{|\alpha|\le
n+2}\Bigl(\|Z^\alpha w_k'(t,\cdot)\|_{L^2(\Omega)}+
\|(1+r)^{-(n-1)/4}Z^\alpha
w_k'\|_{L^2_{s,x}([0,t]\times\Omega)}\Bigr).$$

\par
Our first goal is to inductively prove that if
$\varepsilon<\varepsilon_0$ is sufficiently small, then
\begin{equation}\label{goalM}
M_{k}(T)\le 4C_0\varepsilon
\end{equation}
for every $k=1,2,3,...$.  When $k=1$, \eqref{goalM} follows from
Gronwall's inequality.  We, now, assume that the bound
\eqref{goalM} holds for $k-1$. By Theorem \ref{mainExterior} and
\eqref{reduction}, we then have
\begin{multline}\label{Mk}
M_k(T) \le C\sum_{|\alpha|\le n+2}\int_0^T
\|Z^{\alpha}(1-\eta)(s)Q((u_0+w_{k-1})')(s,\cdot)\|_{L^2(\Omega)}\:ds\\
+C\sup_{0\le s\le T}\sum_{|\alpha|\le n+1}\|Z^\alpha
(1-\eta)(s)Q((u_0+w_{k-1})')(s,\cdot)\|_{L^2(\Omega)}\\
+C\sum_{|\alpha|\le n+1}\|Z^\alpha
(1-\eta)Q((u_0+w_{k-1})')\|_{L^2_{s,x}([0,T]\times\Omega)}\\
+2C\varepsilon + C\sum_{|\alpha|\le
n+2}\int_0^1\|Z^{\alpha}w_k'(s,\cdot)\|_{L^2(\Omega)}\:ds.
\end{multline}
Let's examine the pieces on the right separately.
\par
Since $Q$ is quadratic, for $|\alpha|\le n+2$, we have
\begin{multline}\label{Qquad}
|Z^{\alpha}Q((u_0+w_{k-1})')(s,x)|\\ \le C\left(\sum_{|\alpha|\le
n+2} |Z^\alpha
(u_0+w_{k-1})'(s,x)|\right)\left(\sum_{|\alpha|\le\frac{n+2}{2}}
|Z^\alpha (u_0+w_{k-1})'(s,x)|\right).
\end{multline}

Thus, by \eqref{inftyEstimate} and the standard Sobolev lemma, for
$j=1,2,3,...$, we have
\begin{multline*}
\|Z^{\alpha}Q((u_0+w_{k-1})')(s,x)\|_{L^2(\{2^j\le |x|\le
2^{j+1}\})}\\
\le C\sum_{\alpha\le n+2}\|r^{-(n-1)/4} Z^\alpha
u'_0(s,\cdot)\|^2_{L^2(\{2^{j-1}\le |x|\le 2^{j+2}\})}\\
+C\left(\sum_{|\alpha|\le n+2}\|Z^\alpha
u'_0(s,\cdot)\|_{L^2(\Omega)}\right)\left(\sum_{|\alpha|\le
n+2}\|Z^{\alpha}w'_{k-1}(s,\cdot)\|_{L^2(\{2^j\le |x|\le 2^{j+1}\})}\right)\\
+C\left(\sum_{|\alpha|\le n+2}\|Z^\alpha
u'_0(s,\cdot)\|_{L^2(\{2^j\le |x|\le
2^{j+1}\})}\right)\left(\sum_{|\alpha|\le
n+2}\|Z^{\alpha}w'_{k-1}(s,\cdot)\|_{L^2(\Omega)}\right)\\
 +C\sum_{|\alpha|\le
n+2}\|r^{-(n-1)/4}Z^{\alpha}w'_{k-1}(s,\cdot)\|^2_{L^2(\{2^{j-1}\le
|x|\le 2^{j+2}\})}.
\end{multline*}
Since $u_0(s,x)$ vanishes for $s>1$, applying \eqref{reduction},
gives
\begin{equation}\label{piece1}
\int_0^T
\|Z^{\alpha}(1-\eta)(s)Q((u_0+w_{k-1})')(s,\cdot)\|_{L^2(\Omega)}\:ds
\le C(C_0\varepsilon+M_{k-1}(T))^2.
\end{equation}
\par
For the second term on the right of \eqref{Mk}, by \eqref{Qquad}
and the standard Sobolev lemma, for $|\alpha|\le n+1$, we have
\begin{multline*}
\|Z^\alpha (1-\eta)(s)Q((u_0+w_{k-1})'(s,\cdot)\|_{L^2(\Omega)}\le
C \sum_{n+1} \|Z^{\alpha} u_0'(s,\cdot)\|^2_{L^2(\Omega)}\\
+2C\left(\sum_{|\alpha|\le n+1}\|Z^\alpha
u'_0(s,\cdot)\|_{L^2(\Omega)}\right)\left(\sum_{|\alpha|\le
n+1}\|Z^{\alpha}w'_{k-1}(s,\cdot)\|_{L^2(\Omega)}\right)\\
 +C\sum_{|\alpha|\le
n+1}\|Z^{\alpha}w'_{k-1}(s,\cdot)\|^2_{L^2(\Omega)}.
\end{multline*}

Thus, by \eqref{reduction}, we have
\begin{equation}\label{piece2}
\|Z^\alpha (1-\eta)(s)Q((u_0+w_{k-1})'(s,\cdot)\|_{L^2(\Omega)}
\le C (C_0\varepsilon + M_{k-1}(T))^2.
\end{equation}

\par Finally, for the third term on the right side of \eqref{Mk},
again by \eqref{Qquad},
\begin{multline*}
\|Z^{\alpha}Q((u_0+w_{k-1})')(s,x)\|_{L^2(\{2^j\le |x|\le
2^{j+1}\})}\\
\le C\sum_{\alpha\le n+2}\|r^{-(n-1)/4} Z^\alpha
u'_0(s,\cdot)\|^2_{L^2(\{2^{j-1}\le |x|\le 2^{j+2}\})}\\
+C\left(\sum_{|\alpha|\le n+2}\|Z^\alpha
u'_0(s,\cdot)\|_{L^2(\Omega)}\right)\left(\sum_{|\alpha|\le
n+2}\|Z^{\alpha}w'_{k-1}(s,\cdot)\|_{L^2(\{2^j\le |x|\le 2^{j+1}\})}\right)\\
+C\left(\sum_{|\alpha|\le n+2}\|Z^\alpha
u'_0(s,\cdot)\|_{L^2(\{2^j\le |x|\le
2^{j+1}\})}\right)\left(\sum_{|\alpha|\le
n+2}\|Z^{\alpha}w'_{k-1}(s,\cdot)\|_{L^2(\Omega)}\right)\\
 +C\sum_{|\alpha|\le
n+2} 2^{-j(n-1)/4}
\|r^{-(n-1)/4}Z^{\alpha}w'_{k-1}(s,\cdot)\|_{L^2(\{2^{j-1}\le
|x|\le 2^{j+2}\})}\|Z^{\alpha}w'_{k-1}(s,\cdot)\|_{L^2(\Omega)}.
\end{multline*}

Thus,
\begin{equation}\label{piece3}
\|Z^\alpha
(1-\eta)Q((u_0+w_{k-1})'\|_{L^2_{s,x}([0,T]\times\Omega)} \le C
(C_0\varepsilon + M_{k-1}(T))^2.
\end{equation}

\par
By combining \eqref{Mk}, \eqref{piece1}, \eqref{piece2}, and
\eqref{piece3}, we see that
$$M_k(T)\le 3C(C_0\varepsilon+M_{k-1}(T))^2+2C\varepsilon + C
\sum_{|\alpha|\le n+2}\int_0^1 \|Z^\alpha
w'_k(s,\cdot)\|_{L^2(\Omega)}\:ds.$$ Thus, if $\varepsilon$ is
small enough, \eqref{goalM} follows from Gronwall's Inequality.

\par
Furthermore, if we set
\begin{multline*}
A_k(T)=\sup_{0\le t\le T}\left(\sum_{|\alpha|\le n+2}\|Z^\alpha
(w'_k-w'_{k-1})(t,\cdot)\|_{L^2(\Omega)}\right)\\
+\sum_{|\alpha|\le n+2}\|(1+r)^{-(n-1)/4}Z^\alpha
(w'_k-w'_{k-1})\|_{L^2_{s,x}([0,T]\times\Omega)}
\end{multline*}
and argue as in Section 3, we see that
$$A_k(T)\le \frac{1}{2}A_{k-1}(T)$$
if $\varepsilon$ is small enough.
\par
We have, thus, shown that $w_k$ converge to a solution of
\eqref{wEquation} which satisfies
$$\sum_{|\alpha|\le n+2} \left(\|Z^\alpha
w'(t,\cdot)\|_{L^2(\Omega)} + \|(1+r)^{-(n-1)/4}Z^\alpha
w'\|_{L^2_{s,x}([0,t]\times\Omega)}\right)\le C\varepsilon$$ for
any $t$.  Thus, $u=u_0+w$ is a solution to \eqref{main} satisfying
an analogous bound.  If the data satisfies the compatibility
conditions to infinite order, the solution will be $C^\infty$ on
$\Rplus\times\Omega$ by standard local existence results (see
e.g., \cite{KSS1}).  This completes the proof of Theorem
\ref{GEexterior}.\qed


\end{document}